%
%
%
%
%
%
%

\input amstex
\input epsf
\documentstyle{amsppt}
\nologo
\TagsOnRight
\magnification 1200
\voffset -1cm
\def\epsfsize#1#2{\hsize}

\def\on#1{\expandafter\def\csname#1\endcsname{{\operatorname{#1}}}} \on{im}
\on{lf} \on{id} \on{dist} \on{deg} \on{Hom}  \on{holink}
\on{eq} \on{Tel} \on{Cyl} \on{Fr}  \on{Int} \on{Bl}
\def\fn#1{\expandafter\def\csname#1\endcsname{\operatorname{#1}\def\shortcut
{\if(\next{}\else\,\fi}\futurelet\next\shortcut}} \fn{coker} \fn{mod} \fn{rel}
\let\emb\hookrightarrow \let\imm\looparrowright  
\def\R{\Bbb R} \def\?{@!@!@!@!@!}  \def\z#1#2{#1\?#2\?}
\def\Z{{\Bbb Z}\def\shortcut{\if/\next{}\z\fi}\futurelet\next\shortcut}
\let\tl\tilde \let\but\setminus \let\x\times \let\eps\varepsilon \on{coind}
 \let\phi\varphi \def\T{_{\?\text{\fiverm T}}}

 \on{rk} \fn{lk} \fn{st}
  \on{mesh} \on{fr}
  
\def\Y{_{\sssize\therefore}} \def\y{_{\sssize\Delta}}
 \let\iff\Leftrightarrow
\def\U#1{^{\left<#1\right>}}

\topmatter 
\thanks Partially supported by the Russian Foundation for Basic Research Grant 
No.\ 15-01-06302
\endthanks
\address Steklov Mathematical Institute of Russian Academy of Sciences,
8 Gubkina St., Moscow 119991, Russia \endaddress
\email melikhov\@mi.ras.ru \endemail

\title Gauss-type formulas for link map invariants \endtitle

\author Sergey A. Melikhov \endauthor


\abstract\nofrills
We find that Koschorke's
$\beta$-invariant and the triple $\mu$-invariant of link maps in the critical
dimension can be computed as degrees of certain maps of configuration
spaces --- just like the linking number.
Both formulas admit geometric interpretations in terms of Vassiliev's
ornaments via new operations akin to the Jin suspension, and both were
unexpected for the author, because the only known direct ways to extract
$\mu$ and $\beta$ from invariants of maps between configuration spaces
involved some homotopy theory (Whitehead products and the stable Hopf
invariant, respectively).
\endabstract
\endtopmatter

\document
\head Introduction \endhead

The linking number of a link $S^1\sqcup S^1\subset\R^3$ is the degree of
the Gauss map $S^1\x S^1\to\R^3\x\R^3\but\Delta\simeq S^2$.
A considerable progress has been made in finding similar expressions for
other invariants of knots and links and higher-dimensional embeddings, and
the question now arises, do there ever exist two embeddings of the same
manifold into $\R^m$ that cannot be distinguished by a so expressible
invariant?
We shall now see that this question, to be made more precise in a second,
is a natural generalization of the Vassiliev Conjecture on completeness of
finite type invariants of classical knots.

\definition{Configuration space integrals}
It is known that all rational Vassiliev invariants of knots can be computed
from ``configuration space integrals'' (see \cite{Vo2} and references there).
These are roughly the degrees of the maps between the Fulton--MacPherson
compactification $(S^1)^{[r]}$ (see definition at the end of the
introduction) of the configuration space $(S^1)^{(r)}$ of $r$-tuples of
distinct points of $S^1$ and manifolds homotopy equivalent to $(\R^3)^{(r)}$
with some diagonals filled in, for $r=2,3\dots$.
More precisely, since each $(S^1)^{[r]}$ is a manifold with nonempty
boundary, the integrals are only well-defined when summed up with
``correction terms'', which involve ``joint'' configuration spaces, where
some points are constrained to be on the knot yet are not allowed to collide
with additional points running over the $3$-space.
Necessity of correction terms is readily clear from the fact that the
Fulton--MacPherson compactification is a PL invariant (see its definition
below) --- but every PL knot in $\R^3$ is PL isotopic (i.e\.
non-locally-flat PL homotopic through embeddings) to the unknot.

Another approach to extracting Vassiliev invariants from
configuration spaces is found in \cite{BCSS}.
It follows from results of D. Sinha \cite{Si} and I. Voli\' c \cite{Vo1}
that all rational type $n$ Vassiliev invariants of a string knot
$(I,\partial)\emb(I^3,\partial)$ are contained in the {\it aligned}
$S_{2n}$-equivariant stratum preserving homotopy class of the induced map
between the Fulton--MacPherson compactifications $I^{[2n]}$ and
$(I^3)^{[2n]}$.
(Technically, here $(I^3)^{[2n]}$ must be endowed with tangent vectors, but
this can be remedied by increasing $2n$ to $4n$.)
As indicated above, this result couldn't have been true without the word
``aligned'', so let us recall its meaning from \cite{Si}.
A point $x$ of $(I^3)^{[r]}$ represents an $r$-tuple $(x_1,\dots,x_r)$ of
(not necessarily distinct) points of $I^3$ along with some additional data
for each substring $s_j$ of coincident points $x_{i_1}=\dots=x_{i_{k_j}}$.
Thinking of $x$ as the limit of a curve of points $x(t)$, $t\to\infty$, in
the configuration space, these data can be derived from the relative
velocities of the moving points $x_{i_1}(t),\dots,x_{i_{k_j}}(t)$
(the details can be found in \cite{FM}).
If $v_{i_k}$ is the tangent vector to the curve $x_{i_k}(t)$ at $t=\infty$
(which must exist for $x$ to be defined), $x$ is called aligned if
$v_{i_1},\dots,v_{i_{k_j}}$ are all collinear for each $j$.
In particular, all points in $(I^3)^{(r)}$ are aligned since each $k_j=1$.
A map $I^{[r]}\to (I^3)^{[r]}$ is termed aligned if all points in its image
are aligned.
Every map $I^{[r]}\to (I^3)^{[r]}$ yielded by a genuine knot $I\emb I^3$ is
aligned simply because $I$ is one-dimensional.

On the other hand, there is a direct geometric argument due to J. Conant
showing that knots that are $C_n$-equivalent ($\iff$ indistinguishable by
invariants of type $<n$) yield maps $I^{[n]}\to (I^3)^{[n]}$ that are
connected by an aligned $S_n$-equivariant stratum preserving homotopy.
It also works to show that this homotopy $I^{[r]}\to(I^3)^{[r]}$ can be
chosen to agree with analogous homotopies $I^{[i]}\to(I^3)^{[i]}$,
(a precise definition of such agreement is given at the end of
the introduction).
\enddefinition

\definition{Polynomial compactification}
Motivated by entirely unrelated considerations (resolution of singularities
of smooth maps), the present author has defined a new compactification
$M^{\{r\}}$ of the configuration space $M^{(r)}$ of a compact smooth manifold
$M$, which can be viewed as a polynomial analogue of the ``linear''
Fulton--MacPherson compactification, and which is obtained by successfully
blowing it up along submanifolds of increasingly degenerate configurations,
as measured by coranks of the collections of vectors
$v_{i_1},\dots,v_{i_{k_j}}$ \cite{M*}.
In particular, aligned maps of \cite{BCSS} are precisely those
($S_r$-equivariant stratum preserving) maps between the Fulton--MacPherson
compactifications that lift to the polynomial compactifications.

Thus, \#1 below in the case of classical knots is nearly equivalent to
(follows from the rational version of and implies) the Vassiliev Conjecture.
\enddefinition

\proclaim{Conjecture} If $m-n\ge 3$, let $r$ be the maximal number such that
$m>\frac{(r+1)(n+1)}r$, i.e\. there are no $r$-tuple points in a generic
$r$-parametric homotopy $I^n\to I^m$.
If $m-n=2$, let $r$ stand for an arbitrarily large finite number depending
on the given embeddings.

1. Smooth embeddings between compact smooth manifolds $N^n\emb M^m$ are
classified up to smooth isotopy by homotopy of the induced maps
$N^{\{r\}}\to M^{\{r\}}$ through $S_r$-equivariant stratum preserving
maps agreeing with those at all levels $r'<r$.

2. PL embeddings $X^n\emb M^m$ of a compact polyhedron into a compact PL
manifold are classified up to PL isotopy by homotopy of the induced maps
$N^{[r]}\to M^{[r]}$ through $S_r$-equivariant stratum preserving
maps agreeing with those at all levels $r'<r$.

3. Topological embeddings between compact manifolds $N^n\emb M^m$ are
classified up to topological isotopy by homotopy of the power maps
$N^r\to M^r$ through $S_r$-equivariant maps whose restriction to each
diagonal $\Delta^S_N$, $S\i\{1,\dots,r\}$ is $S_{|S|}$-equivariant and
sending the complement to every $\Delta^S_N$ into the complement to
$\Delta^S_M$.

4. Link maps $S^{n(1)}\sqcup\dots\sqcup S^{n(k)}\to\R^m$ are classified
up to link homotopy by homotopy classes of the induced maps
$S^{n(i_1)}\x\dots\x S^{n(i_s)}\to (\R^m)^s\but\bigcup_{p,q}\Delta^{i_p=i_q}$,
where $s\le r$ and $\Delta^{i=i}$ stands for $\emptyset$.
\endproclaim

By the results of Haefliger and Weber--Harris \cite{Har} the answer is
certainly yes to \#1 and \#2 in the metastable range $m>\frac{3(n+1)}2$,
$n>1$.
By results of Habegger--Kaiser and Koschorke/Massey, \#4 is known to hold
for $r\le 3$.
Validity of \#3 for non-triangulable manifolds in $r=2$ follows from
\cite{MS}, where the Haefliger--Weber criterion is generalized to arbitrary
compacta.

A homotopy theoretic classification of embeddings with $r=3$ was given in
\cite{GKW} using Goodwillie's Calculus of Functors.
It may possibly follow from this (cf\. \cite{BCSS}) that the answer to \#1
is affirmative for $r=3$ as well.
However, since the classification of \cite{GKW} involves taking limits of
diagrams indexed by all (uncountably many) points of the manifold, it
should be more valuable, from a geometric standpoint, if \#1 is proved
directly rather than deduced from such results.

As for \#2, by an argument from \cite{M3} it implies the ``links modulo
knots'' version of the Vassiliev Conjecture \cite{M3}.

In codimension two, there is a natural weaker version of \#1 and \#2,
involving simultaneous homotopies on all levels $r=1,2,\dots$, agreeing
for different levels, or, equivalently, a homotopy between the
infinite telescopes.
The difference between this and the original version may be captured by
a $\lim^1$ obstruction.
A similar weaker version can be stated for \#3.

For embeddings into $\R^m$ rather than other manifolds, the thinnest
diagonal $\Delta^{1=\dots=r}$ in \#3 and the corresponding strata of
the compactifications in \#1 and \#2 can be deleted both in the domain
and range.
Other diagonals cannot be deleted:

\definition{Quasi-isotopy and Whitehead link}
Every embedding obstruction from the deleted cube in fact obstructs
quasi-embeddability, i.e\. existence for each $\eps>0$ of a map
$g_\eps\:X\to\R^m$ with point-inverses of diameter $<\eps$ (such a map is
called an $\eps$-map).
Indeed, if $\eps$ is sufficiently small, $g_\eps$ yields an
$S_3$-equivariant map between the simplicial deleted cubes, which are
$S_3$-equivariantly homotopy equivalent to the ordinary deleted cubes.
Similarly, every isotopy invariant from the deleted cube is actually
invariant under quasi-isotopy, i.e\. existence for each $\eps>0$ of
a homotopy through $\eps$-maps between given the embeddings.

Haefliger's Whitehead link $W\:S^{2k-1}\sqcup S^{2k-1}\emb\R^{3k}$ can be
obtained from the Borromean rings
$B\:S^{2k-1}\sqcup S^{2k-1}\sqcup S^{2k-1}\emb\R^{3k}$ (see \S2)
by connecting any two components by a thin tube.
Assuming $k>1$, it does not matter which tube is chosen, since the complement
is simply-connected.
Since $B$ is Brunnian, the tubed component is null-homotopic in
the complement to the unmodified one.
Hence $B$ is link homotopic and so quasi-isotopic to the unlink.
On the other hand, $W$ is not isotopic to the unlink if $n\ne 1,3,7$.
Indeed, the homotopy class of the untubed component in the complement to
the tubed component is the Whitehead product
$[\iota,\iota]\in\pi_{2k-1}(S^k)$ (see \cite{Hae}), which is nonzero
when $k\ne 1,3,7$.
So $W$ is not ambient isotopic to the unlink for such $k$, hence not PL
isotopic to the unlink by Zeeman's theorem and hence not topologically
isotopic to the unlink by Edwards' theorem (see references e.g\. in
\cite{M2}).
See \cite{Sk$_{\text A}$} for further examples of this type.
\enddefinition

The version of \#3 omitting the condition that the restrictions to
the diagonals be equivariant was stated in \cite{M2}.
However, this condition, as well the conditions of agreement with levels
$<r$ in \#1 and \#2, turn out to be necessary for $r=3$ by the following
example.

\definition{Semi-contractible Whitehead links}
In the cases $n=3,7$ each component of the Whitehead link
$W\:S^{2k-1}\sqcup S^{2k-1}\emb\R^{3k}$ is null-homotopic in the complement
to the other one (see above).
Combining the tracks of such two null-homotopies yields a link map
$S^{2k}\sqcup S^{2k}\to\R^{3k+1}$ with a non-zero $\beta$-invariant.
Since the $\beta$-invariant vanishes on such link maps with one component
embedded (see \S3), it does not depend on the choice of the homotopies
and so is an isotopy invariant of $W$.
Thus $W$ is non-trivial for $k=3,7$ as well.

Now to get an isovariant homotopy between the power maps
$(S^{2k-1}_1\sqcup S^{2k-1}_2)^3\to(\R^{3k})^3$ of the Whitehead link and
a trivial link we only need to construct it separately for each connected
component of the domain.
It is easy to define it on $(S^{2k-1}_1)^3$ and $(S^{2k-1}_2)^3$.
Now using a null-homotopy of $S^{2k-1}_1$ in the complement to $S^{2k-1}_2$,
we can also define it on $S^{2k-1}_1\x (S^{2k-2}_2)^2$, and vice versa.
\enddefinition

\comment

The results of \S2 and \S3 open up a way to a possibly
complete elimination of homotopy theory (such as Whitehead products and
in particular group commutators) from the theory of configuration space
invariants of link maps in favor of cohomology (or generalized cohomology
when appropriate, cf\. \S4).
This reduction may be not without a little pain (as in the proof of
Lemma 2.4) which hopefully can be converted into pleasure by someone
knowledgeable of operads or $n$-categories.
Perhaps our approach can be understood as a configuration space analogue
of computing Massey products on the link complement.
\endcomment

The possibility of coincidence of indices in \#4 is essential: without
allowing it, we would be dealing with Koschorke's $\kappa$-invariant, that is,
the homotopy class of the map $S^{n(1)}\x\dots\x S^{n(k)}\to(\R^m)^{(k)}$.
But for link maps $S^2\sqcup S^2\to\R^4$ it is strictly weaker than Kirk's
$\sigma$-invariant.
Correspondingly, a notable special case of \#4 concerns Koschorke's
$\tl\beta$-invariant, and in particular Kirk's $\sigma$-invariant of link
maps $S^2\sqcup S^2\to\R^4$.

\remark{Remark} In the last paragraph of the paper \cite{Ko4}, Koschorke
mentions a ``possibly new additive invariant'' of link maps
$S^2\sqcup S^2\to\R^4$ taking values in ``a cyclic group'', which ``might
detect nontrivial elements in the kernel of $\sigma$'' (the possible source
of this invariant can be inferred from \cite{Ko4}).
The present author has shown in \cite{MR; Theorem 2.6(b)}, but failed to
state explicitly there, that the ``cyclic group'' in question is trivial
after all.
\endremark

\definition{Outline of the paper and discussion}
Problems 1--4 are concerned with a homotopy-theoretic classification
of embeddings and link maps.
In view of the Bott--Taubes integrals, one may also wonder
whether a classification in terms of (generalized equivariant) cohomology
invariants of maps between configuration spaces is possible.
As observed in \cite{M5}, this is indeed the case when $r=2$ (for smooth and
PL embeddings and link maps).
In this paper it is proved that such a classification is also possible for
link maps with $r=3$, by expressing the triple $\mu$-invariant and
the $\beta$-invariant of link maps as (generalized equivariant) degrees of
certain maps between configuration spaces.
The only previously known ways of extracting $\mu$ and $\beta$ from
configuration spaces involved a Whitehead product (in the case of $\mu$)
and the stable Hopf invariant (in the case of $\beta$), even in those
dimensions where our formulae involve only ordinary (equivariant) cohomology.

These formulas for $\mu$ and $\beta$ are geometrically well understood.
Geometrically, the link map is converted in both cases into Vassiliev's
``ornament'', which is a map $S^i\sqcup S^j\sqcup S^k\to\R^m$ with no
triple points involving all three components (any intersections and
self-intersections of any pair of components are allowed).
In the range $r=3$, ornaments are classified (up to homotopy through
ornaments in the case of $\mu$, and up to homotopy through ``equivariant
ornaments'' in the case of $\beta$) by a cohomological invariant directly
analogous to the linking number.

These two geometric constructions producing ornaments from link maps are
in fact simplest instances of a general method, allowing to ``abelianize''
a lower-dimensional, ``non-linear'' situation into a higher-dimensional
``linear'' one.
A third instance, producing link maps from semi-contractible links, was
known as ``Jin suspension'' since mid 80s, and is used e.g\. in the
author's new proof of Kirk--Livingston and Nakanishi--Ohyama theorems
\cite{M6}.
It is expected that this method, further instances of which are now known,
will be a promising substitute for the Whitney tower/grope/$C_k$-move
techniques, which nobody managed to make work adequately in
the multi-component case as yet.
\enddefinition

\definition{Fulton--MacPherson compactification}
Let $X$ be any compact polyhedron.
What we call the {\it Fulton--MacPherson compactification} $X\U2$
of the two-point configuration space $X^{(2)}=X\x X\but\Delta$ is
PL homeomorphic to $X\x X\but\Int N$, where $N$ is an equivariant regular
neighborhood of $\Delta$.
Using an equivariant collapse of $N$ onto $\Delta$, one can define a map
$\pi\:X\U2\to X\x X$ such that $\pi^{-1}(\Delta)=\Fr N$ and $\pi$ sends
$\Int X\U2$ homeomorphically onto $X^{(2)}$ --- thus indeed making $X\U2$
a compactification of $X^{(2)}$.
In contrast to the usual definition of $X\U2$ in the smooth case, our
definition does not uniquely describe $X\U2$ up to a homeomorphism commuting
with $\pi$, nor even up to a homeomorphism fixed on $X^{(2)}$.
But it does uniquely define $\pi$ up to PL ``block equivalence'', which is
a notion from the same category as PL block bundle isomorphism.
That is, given an equivariant triangulation of $X\x X$, any two instances
$X\U2_\alpha$, $X\U2_\beta$ of $X\U2$ are related by an equivariant PL
homeomorphism taking the preimage in $X\U2_\alpha$ of each dual cone of
the triangulation of $X\x X$ into its preimage in $X\U2_\beta$; moreover
the homeomorphism strictly commutes with $\pi$ over dual cones disjoint
from $\Delta$.

For $r=2,3$ the Fulton--MacPherson compactification of the $r$-point
configuration space $X^{(r)}$ coincides with its simplified version
$X\U{r}$, defined as follows.
Let us consider an $S_r$-regular neighborhood $N_1$ of the thinnest diagonal
$\Delta^{1=\dots=r}$ of $X^r=X\x\dots\x X$; an $S_r$-equivariant regular
neighborhood $N_2$ of union of the next-thick diagonals
$\Delta^{i_1=\dots=i_{r-1}}$ and $\Delta^{i_1=\dots=i_j;\,i_{j+1}=\dots=i_r}$
relative to $N_1$; and so on, ending up with an $S_r$-equivariant regular
neighborhood $N_{r-1}$ of union of the thickest diagonals $\Delta^{i=j}$
relative to $N_1\cup\dots\cup N_{r-2}$.
(The diagonals are indexed by all partitions of $\{1,\dots,r\}$, except for
the partition into singletons.
The thickness of a diagonal is the number of blocks in the partition.)
Now $X\U{r}$ is PL homeomorphic to $X^r\but\Int (N_1\cup\dots\cup N_{r-1})$,
and is endowed with the collection of codimension zero subpolyhedra
$\partial_i X\U{r}=\Fr N_i\cap\partial X\U{r}$ of the corona
$\partial X\U{r}=\Fr (N_1\cup\dots\cup N_{r-1})$.
Connected (normally) components $\partial_\Pi X\U{r}$ of each
$\partial_i X\U{r}$ are indexed by partitions $\Pi$ with $|\Pi|=i$ blocks.
If $X$ is a closed manifold, $X\U{r}$ is a manifold with boundary
$\partial X\U{r}$ and with corners
$\partial_{\Pi_1} X\U{r}\cap\dots\cap\partial_{\Pi_k} X\U{r}$, where
each $\Pi_i$ refines $\Pi_{i+1}$.
The projection $\pi\:X\U{r}\to X^r$ is well-defined up to block equivalence,
and its restriction to each $\partial_\Pi X\U{r}$ factors through a map
$\pi_\Pi\:\partial_\Pi X\U{r}\to X\U{|\Pi|}$, well-defined up to block
equivalence, and the projection $\pi\:X\U{|\Pi|}\to X^{|\Pi|}$.
Gluing together $X\U{r}$ and the cylinders of all $\pi_\Pi$'s, we recover
the product $X^r$.

Restricting our attention to the diagonals indexed by subsets (which can be
identified with partitions into one subset and singletons) in the above
definition, we get the {\it Fulton--MacPherson compactification} $X^{[r]}$.
The role of $\partial_\Pi X\U{r}$'s is played by the images
$\partial_S X^{[r]}$ under the natural projection $X\U{r}\to X^{[r]}$ of
$\partial_{\{S\}} X\U{r}$ for all subsets $S\i\{1,\dots,r\}$.
A map $f\:X^{[r]}\to Y^{[r]}$ will be called {\it stratum preserving} if
$f(\partial_S X^{[r]})=\partial_S Y^{[r]}$ for each subset $S$.
(It would perhaps be more accurate to call $f$ preserving closures of strata,
but ``stratum preserving'' appears to be a standard terminology \cite{Si}.)
One can define the maps $\pi_S\:\partial_S X^{[r]}\to X^{[r+1-|S|]}$
similarly to the above.
We shall say that an $S_r$-equivariant stratum preserving map
$f_r\:X^{[r]}\to Y^{[r]}$ {\it agrees with similar lower level maps} if
there exist $S_i$-equivariant stratum preserving maps
$f_i\:X^{[i]}\to Y^{[i]}$ for all $i<r$ such that
$\pi_S f_i=f_j\pi_S$ for all subsets $S$ of $\{1,\dots,i\}$ with $j=i+1-|S|$.

If $X$ is a smooth manifold or $r=2$, each $\partial_S X^{[r]}$ and
$\partial_i X\U{r}$ admits a canonical involution, and these can be used
to make $X^{[r]}$ and $X\U{r}$ into closed manifolds --- the original
projective versions of the Fulton--MacPherson and the
Fulton--MacPherson--Ulyanov compactifications \cite{FM}, \cite{U}.
\enddefinition

\bigskip

I am grateful to P. Akhmetiev, V. M. Nezhinskij for useful discussions.
This paper arose in connection with renewal of my interest in extracting link
invariants from configuration spaces, for which I should thank participants
of the Workshop `Moduli spaces of knots' (Palo Alto, January 2006).

\bigskip

\definition{Note added in 2017} This preprint was privately circulated and presented at conferences and seminars by the autor in 2006-07 (Alexandroff Readings, Moscow, 
May 2006; UF-FSU Topology Conference, Gainesville, December 2006; Dartmouth College
topology seminar, April 2007; etc.)
I hesitated to publish it at that time as I hoped to get more progress on 
the conjectures stated in the introduction; but numerous other projects have been
distracting me from this task so far.

In recent years, a construction pretty similar to (but not exactly same as)
the triple-point Whitney trick in Theorem 1.1 below has been independently discovered
by I. Mabillard and U. Wagner, who successfully employed it to obtain nice results.
(They call theirs the ``triple Whitney trick'', but I prefer to reserve this title for
a more elaborate construction, involving the triple-point Whitney trick as one of 
several steps; it can be used, in particular, to obtain a geometric proof of
the Habegger--Kaiser classification of link maps in the 3/4 range --- as presented in 
my talk at the Postnikov Memorial Conference in Bedlewo, June 2007, --- and will hopefully
appear elsewhere.)
As the triple-point Whitney trick has attracted considerable interest, I have added 
several footnotes clarifying the proof of Theorem 1.1 without altering its original
(admittedly, somewhat sketchy) text from 2006.
\enddefinition

\definition{Notation}
The following notation will be used throughout the paper.
Let $X$ be a compact $n$-polyhedron.
$\bar X$ will denote the quotient of the {\it deleted product}
$\tl X:=X\x X\but\Delta_X$ by $\Z/2$, acting by the factor exchanging
involution.

If $X$ is a polyhedron, let $\tl X\Y$ denote the {\it deleted cube}
$X\x X\x X\but(\Delta^{1=2}\cup\Delta^{2=3}\cup\Delta^{3=1})$, where each
$\Delta^{i=j}$ denotes the thick diagonal $\{(x_1,x_2,x_3)\mid x_i=x_j\}$.
Also let $\tl X\y$ denote $X\x X\x X\but\Delta^{1=2=3}$, where
$\Delta^{1=2=3}$ denotes the thin diagonal $\{(x,x,x)\mid x\in X\}$.
\enddefinition

\head 1. Ornaments ($2k-1$ in $3k-1$) \endhead

Doodles were originally defined by Fenn and Taylor \cite{FT} (see also
\cite{Fe}) as triple point free maps of a collection of circles into
the plane that embed each circle.
In our terminology we follow Khovanov \cite{Kh}, who redefined them by
dropping the condition that each component be embedded.
It has been established that Khovanov's doodles are classified up to
{\it doodle homotopy}, i.e\. homotopy through doodles, by their finite type
invariants \cite{Mer2}.

\definition{Invariant of ornaments}
If $f=f_1\sqcup f_2\sqcup f_3\:X_1\sqcup X_2\sqcup X_3\to\R^m$ is an
{\it ornament}%
\footnote{Ornaments were introduced by Vassiliev \cite{Va1} as a modification
of Fenn and Taylor's doodles \cite{FT}; see also \cite{Mer1} (especially \S6
and Appendix D) and \cite{Va2; Chapter 6}.}%
, that is $f(X_1)\cap f(X_2)\cap f(X_3)=\emptyset$ (but $f$
may have triple points), let $$\breve\mu(f)\in H^{2m-1}(X_1\x X_2\x X_3)$$ be
the image of a fixed generator of $H^{2m-1}(S^{2m-1})$ under the composition
$$X_1\x X_2\x X_3@>f_1\x f_2\x f_3>>\tl\R^m\y@>\simeq>>S^{2m-1}.$$
Then $\breve\mu(f)$ is invariant under {\it ornament homotopy} (i.e\.
homotopy through ornaments).
We note that $\breve\mu$ depends on the ordering of $X_i$.
Since $S_3$ acts on $H^{2m-1}(\tl\R^m\y)$ by the sign homomorphism,
$\breve\mu$ remains unchanged when $X_1$, $X_2$, $X_3$ are permuted cyclicly,
and reverses the sign when any two are interchanged.

If each $X_i$ is an oriented connected closed $(2k-1)$-manifold and $m=3k-1$,
then $\breve\mu(f)$ is an integer, and it clearly equals the algebraic number
of {\it $1=2=3$ points}, i.e\. triple points of intersection between $X_1$,
$X_2$ and $X_3$, in a generic homotopy from $f$ to $t$.
In particular, in the case of ornaments $f\:S^1\sqcup S^1\sqcup S^1\to\R^2$
it is an extension of the $\mu$-invariant of \cite{Fe}.
Note that in this case it is incomplete up to ornament homotopy (see e.g\.
\cite{Va1}) or ornament concordance, but is complete up to ornament
cobordism (see \cite{Fe}).
\enddefinition

\definition{Borromean ornament}
Let us think of $S^{3k-1}$ as the unit sphere in $\R^{3k}$.
Let $b$ be a triple point free map
$S^{2k-1}\sqcup S^{2k-1}\sqcup S^{2k-1}\to S^{3k-1}$ embedding the components
onto $S^{3k-1}\cap\R^k\x\R^k\x 0$, $S^{3k-1}\cap\R^k\x 0\x\R^k$ and
$S^{3k-1}\cap 0\x\R^k\x\R^k$.
(This was called the Borromean doodle in \cite{Fe} when $k=1$.)
The obvious null-homotopy of $b$ has one $1=2=3$ point, so \
$\breve\mu(b)=\pm 1$.
\enddefinition

\proclaim{Theorem 1.1} $\breve\mu(f)$ is a complete invariant of ornament
homotopy if each $X_i$ is an orientable $(2k-1)$-manifold and $m=3k-1$,
$k>2$.
\endproclaim

\demo{Proof}
Let us assume for simplicity that each $X_i$ is connected and closed.
If $\breve\mu(f)=\breve\mu(g)$, there is a homotopy $h\:X\x I\to\R^m\x I$
between the ornaments $f$ and $g$ whose $1=2=3$-points can be paired up with
opposite signs.
Each pair $(p^+,p^-)$ can then be cancelled by a triple Whitney trick.
In more detail, let $p_i^\pm$ be the preimage of $p^\pm$ in $X_i\x I$.
We first arrange that $(p_1^+,p_2^+)$ and $(p_1^-,p_2^-)$ be in the same
component of the double point set $\{(x,y)\mid h(x)=h(y)\}$ in
$X_1\x I\x X_2\x I$ (in case that initially they are not).
To this end we pick points $(q_1^\pm,q_2^\pm)$ in the same components with
$(p_1^\pm,p_2^\pm)$ and such that the double points $f(q_1^+)=f(q_2^+)$ and
$f(q_1^-)=f(q_2^-)$ are not triple points.
Now connect $q_1^+$ and $q_1^-$ by a generic path in $X_1\x I$, disjoint from
preimages of any double points (using that $k>1$) and attach a thin
$1$-handle (i.e\. remove $B^n\x S^0$ and paste in $S^{n-1}\x I$) to
$h(X_2\x I)$ along the image of this path.%
\footnote{In more detail, it may be assumed that the double point set is an oriented manifold that 
is immersed in $X_1\x I$. 
Let us take an oriented connected sum of its components along the original path in $X_1\x I$.
This yields an embedded tube $S^{k-1}\x I$ in $h(X_1\x I)$, and its fiberwise join with the spherical
normal bundle of $h(X_1\x I)$ over the image of the path is the desired embedded tube $S^{2k-1}\x I$
in $\R^m\x I$.}
To restore the topology of $X_2\x I$, we cancel the $1$-handle geometrically
by attaching a $2$-handle along an embedded $2$-disk, meeting $h(X\x I)$
only in the boundary circle (such a disk exists since $k>2$).%
\footnote{Moreover, since $k>2$, the boundary circle may be assumed to be disjoint from $h(X_1\x I)$, so 
that the double point set is unaffected by this step.
Since the normal bundle of the boundary circle in $h(X_2\x I)$ is trivial, it extends over the disk;
the associated sphere bundle is the desired $2$-handle.}
To cancel the $1=2=3$ points, let us connect $(p_1^+,p_2^+)$ and
$(p_1^-,p_2^-)$ by a path%
\footnote{within the double point set}%
, and attach a thin $1$-handle to $h(X_3\x I)$ along
the image of this path.%
\footnote{This $1$-handle is the spherical normal bundle of $h(X_1\x I)\cap h(X_2\x I)$ over the path.
It is attached orientably since the two $1=2=3$ points have opposite signs.}
The topology of $X_3\x I$ can be restored using another $2$-disk like before.%
\footnote{In particular, this $2$-disk is disjoint from $h(X_1\x I\sqcup X_2\x I)$, so no 
new $1=2=3$ points arise.}

Finally, we must either make sure that this construction can be done in a
level-preserving fashion, or alternatively apply the ``ornament concordance
implies ornament homotopy in codimension three'' theorem \cite{M1}. \qed
\enddemo

We note that $\breve\mu(g)=0$ if $g$ is a {\it link map}, i.e\.
$g(X_i)\cap f(X_j)=\emptyset$ whenever $i\ne j$.
Indeed, by first moving $X_1$ away while keeping $f$ fixed on $X_2\cup X_3$,
and subsequently moving $X_2$ away from $X_3$, we get an ornament homotopy
from $g$ to a trivial link map.

This is not a coincidence:

\proclaim{Theorem 1.2} Consider link maps
$f\:S^{n_1}\sqcup\dots\sqcup S^{n_k}\to\R^m$, $n_i<m-1$, and let $\omega(f)$
be the image in the cohomology of $S^{n_1}\x\dots\x S^{n_k}$ of some
cohomology class of $(\R^m)^k\but\bigcup\Delta^{i=j}$.
Suppose that $\omega$ vanishes on any $f$ with one component contained
in an $m$-ball disjoint from the other components.
Then $\omega$ vanishes on all $f$.
\endproclaim

\demo{Proof} Since spherical link maps in codimension two up to link
homotopy (i.e.\ homotopy through link maps) form a group \cite{BT},
\cite{M1}, and $\omega$ is clearly invariant under link homotopy and
additive under connected sum, it suffices to consider the case where $f$
is homotopically Brunnian.
Now by an analysis of Koschorke \cite{Ko1}, \cite{Ko7}, $f_1\x\dots\x f_n$
factors through a Whitehead product, hence is trivial on cohomology. \qed
\enddemo

\definition{$\beta_i$-invariants of animated ornaments}
An {\it animated ornament} is an ornament
$f\:X_1\sqcup X_2\sqcup X_3\to\R^m$ along with three ornament homotopies
starting with $f$, where the $i$th homotopy has support in the $i$th
component and shrinks it to a point.

An animated ornament $(F,f)=\bigsqcup (F_i,f_i)\:\bigsqcup (CX_i,X_i)\to\R^m$
yields a map
$$\Sigma\Y(X_1\x X_2\x X_3)@>\bigsqcup f_i\x f_j\x F_k>>\tl\R^m\y
@>\simeq>>S^{2m-1},$$
which suspends to $S^1*(X_1\x X_2\x X_3)\to\Sigma\Y S^{2m-1}$.
A $\Z[\Z/3]$-module isomorphism $H^{2m}(\Sigma\Y S^{2m-1})\simeq\Z[\zeta]$
can be chosen so that the homotopy between $f_i\x F_j\x f_k$ and
$f_i\x f_j\x F_k$ only meets the support of $\zeta^{i-1}$.
Let us denote the image of $\zeta^{i-1}$ in
$H^{2m}(S^1*(X_1\x X_2\x X_3))\simeq H^{2m-2}(X_1\x X_2\x X_3)$
by
$$\beta_i(F)\in H^{2m-2}(X_1\x X_2\x X_3).$$
By construction, $\beta_i$ is invariant under homotopy through animated
ornaments, remains unchanged if $X_j$ and $X_k$ are interchanged, and
does not depend on $F_i$.
Also $\beta_i+\beta_j+\beta_k=0$ since $1+\zeta+\zeta^2=0$.
\enddefinition

In fact, it follows from the definition that $\beta_i$ is simply the
$\breve\mu$-invariant of the ornament
$\Sigma X_1\sqcup\Sigma X_2\sqcup\Sigma X_3\to\R^{m+1}$
given by $F_j$ and $(f_i\sqcup f_k)\x\id_{[0,1/2]}$ plus some null-homotopies
of $f_i$ and $f_k$ on the upper half of the suspensions, and by
$F_k$ and $(f_i\sqcup f_j)\x\id_{[-1/2,0]}$ plus some null-homotopies
of $f_i$ and $f_j$ on the lower half of the suspensions.

\definition{$\beta_i$-invariants of $\vec v$-animated ornaments}
An ornament $f\:X_1\sqcup X_2\sqcup X_3\to\R^m$ is called
{\it $\vec v$-animated} if no vector of the type $f(x)-q$ is a positive
scalar multiple of $\vec v$, where $q=f(y)=f(z)$ and each $X_i$ contains
precisely one of $x,y,z$.
Clearly, every $\vec v$-animated ornament is animable.

A $\vec v$-animated ornament $f=\bigsqcup f_i\:\bigsqcup X_i\to\R^m$ yields
a map $$X_1\x X_2\x X_3@>\bigsqcup f_i\x f_j\x f_k>>(\R^m)^3\but
\Pi^{\vec v}_1\cup\Pi^{\vec v}_2\cup\Pi^{\vec v}_3@>\simeq>>
\Sigma\Y S^{2m-3}.$$
Then with an appropriately chosen isomorphism
$H^{2m-2}(\Sigma\Y S^{2m-3})\simeq\Z[\zeta]$, the image of $\zeta^{i-1}$ in
$H^{2m-2}(X_1\x X_2\x X_3)$ equals $\beta_i$.
\enddefinition

\definition{$\beta_0$-invariant for $\pm\ne 0$ maps} Let
$f\:X_+\sqcup X_-\sqcup X_0\to\R^{3k-2}$ be
a $\pm\ne 0$ map of a $(2k-2)$-polyhedron (see \S3).
Then any null-homotopy of $f|_{X_1}$, any null-homotopy of $f|_{X_3}$ and
(using the hypothesis on dimensions) any generic null-homotopy of $f|_{X_2}$
combine into an animated ornament.
Moreover, since $\beta_0$ does not depend on the choice of the latter
null-homotopy, it only depends on $f$ and is invariant under $\pm\ne 0$
homotopy of $f$.
Specifically, $\beta_0(f)$ coincides with $\hat\beta(f)$ from \S3 since both
equal the $\breve\mu$-invariant of the ornament $\hat f$ from the proof of
Theorem 3.1.
Alternatively, choosing the animated ornament to be $\vec v$-animated for
some $\vec v$, we see that $\beta_0(f)$ as defined in the previous paragraph
(making use of $\vec v$) equals $\hat\beta(f)$ as defined in \S3 using
the composition $X_1\x X_2\x X_2\to
(\R^m)^3\but(\Pi^{\vec v}_2\cup\Pi^{\vec v}_3)\simeq S^{2m-1}$.
\enddefinition

\head 2. $\mu$-invariant of link maps ($2k-3$ in $3k-3$)  \endhead

The first dimension beyond the metastable range is $n=2k-2$, $m=3k-2$
for embeddability and $n=2k-3$, $m=3k-3$ for isotopy.
For convenience, let us substitute $k-1$ for $k$ and write $n=2k-1$, $m=3k$
from now on.

\definition{Borromean rings}
The Borromean rings $B\:S^{2k-1}\sqcup S^{2k-1}\sqcup S^{2k-1}\emb\R^{3k}$
can be obtained from concentrically embedded spheres
$S^{k-1}_+,S^{k-1}_-\i\R^k$ by forming the joins
$$\multline(S^{k-1}_+\x 0\x 0)*(0\x S^{k-1}_-\x 0)\quad\sqcup\quad
(0\x S^{k-1}_+\x 0)*(0\x 0\x S^{k-1}_-)\\
\quad\sqcup\quad(0\x 0\x S^{k-1}_+)*(S^{k-1}_-\x 0\x 0)\qquad\subset\qquad
\R^k\x\R^k\x\R^k.\endmultline$$
The non-triviality of $B$ cannot be detected by the deleted product, since
every two-component sublink is easily seen to be trivial.
It is well-known, however, that it can be detected by the homotopy class of
one component in the complement to the other two, which is a Whitehead
product \cite{Ha2} and by the triple Massey product of the Alexander duals
to the fundamental homology classes of the components \cite{Ma1}.

As noticed by Koschorke \cite{Ko1}, the non-triviality of $B$ up to link
homotopy can also be detected from the homotopy class of
$F\:S^{2k-1}\x S^{2k-1}\x S^{2k-1}\to\tl\R^{3k}\Y$, which turns
out to factor through the Whitehead product
$S^{6k-1}\to S^{3k-1}\vee S^{3k-1}$.
Indeed, since $B$ restricted to the first two components is link homotopically
trivial, $F$ composed with the projection
$\pi\:\tl\R^{3k}\Y\to\tl\R^{3k}$ that forgets the third point of
the configuration is null-homotopic.
Therefore $F$ is homotopic to a map into a fiber of $\pi$, i.e\. $\R^m$ minus
two points, which collapses onto $S^{3k-1}\vee S^{3k-1}$.

In light of this observation it may look like no triple analogue of
the van Kampen obstruction to isotopy, defined by inducing some cohomology
class of $S_3$, can detect the non-triviality of the Borromean rings
(compare Theorem 1.2).
Nevertheless, we will show in Theorem 2.5 that the triple $\mu$-invariant,
which can be defined using either of the three mentioned approaches to prove
the non-triviality of $B$ (see \cite{Ko7}, \cite{Ko5; 3.11}), in fact equals
a certain $3$-parametric version of the triple van Kampen obstruction.
This agrees with what was apparently expected in \cite{Kr; \S5}.
\enddefinition

\definition{$\mu$-invariant of link maps}
The following definition of the $\mu$-invariant is essentially found in
\cite{Ah; \S4} (in the classical dimension and $\bmod 2$); Akhmetiev
informed the author that he mentioned it as early as 1991 in his talk at
a conference in Kiev.
It may also be viewed as a specific instance of Koschorke's geometric
definition \cite{Ko6}, essentially rediscovered in \cite{MM} (in
the classical dimension), although neither \cite{Ko6} nor \cite{MM}
singles out this specific instance; compare \cite{Ko5; 3.10}.

Let $f\:S^{2k-1}\sqcup S^{2k-1}\sqcup S^{2k-1}\emb\R^{3k}$ be a
{\it Brunnian link}, i.e\. one where each proper sublink is isotopically
trivial.
We define $\mu(f)$ to be the algebraic number of triple points of
intersection between the tracks of a generic null-homotopy of the first
component in the complement to the second; a generic null-homotopy of
the second component in the complement to the third; and a generic
null-homotopy of the third component in the complement to the first.
We will refer to such a triple of null-homotopies as a
{\it $(123)$-null-homotopy} of $f$; also recall that we abbreviate a triple
point of intersection between three distinct components as a
``$1=2=3$ point''.
The three obvious flat disks bounded by the Borromean rings $B$, where
$S^{k-1}_-$ is taken in the bounded complementary domain of $S^{k-1}_+$,
yield a $(123)$-null-homotopy of $B$ with precisely one $1=2=3$ point,
thus $\mu(B)=\pm 1$.
\enddefinition

The following lemma allows to extend the definition to all {\it h-Brunnian}
link maps, i.e\. ones where every proper sub-link-map is link homotopic to
an unlink.

\proclaim{Lemma 2.1} Every h-Brunnian link map $f\:3S^{2k-1}\to\R^{3k}$ is
link homotopic to a link admitting a $(123)$-null-homotopy as well as
a $(321)$-null-homotopy (in fact, to a Brunnian link).
\endproclaim

\demo{Proof} If $k>1$, we first link homotop $f$ to an embedded link using
the Penrose--Whitehead--Zeeman trick (this is the lowest dimension where it
works).
We then make it Brunnian by adding split links resulting from the reflected
$2$-component sublinks of $f$ (using that concordance implies isotopy in
codimension three).

For $k=1$ we could resort to Milnor's classification to get the assertion.
(It might be undesirable to use Milnor's work on $\mu$-invariants in
defining a $\mu$-invariant.
To this end, we can simply unknot the components by a link homotopy, so that
each component becomes null-homotopic in the complement to every other.
This suffices to get the weaker conclusion.) \qed
\enddemo

\proclaim{Theorem 2.2} $\mu(f)$ is well-defined up to link homotopy.
\endproclaim

\demo{Proof} Let $F_1,F_2,F_3\:S^{2k-1}\x I\to\R^{3k}\x I$ be a generic link
homotopy between $f$ and $g$, and let $h_0$ be a $(123)$-null-homotopy of $f$
and $h_1$ a $(321)$-null-homotopy of $g$.
We can extend $h_0$ and $h_1$ to some triple of null-homotopies
$H_1,H_2,H_3\:D^{2k+1}\to\R^{3k}\x I$ of $F_1,F_2,F_3$.
This triple yields an oriented bordism between the $1=2=3$ points of $h_0$
and those of $h_1$ as well as the points of intersection between $H_i$,
$H_j$ and $F_k$ whenever $(ijk)=(123)$.
The latter are the $F_k$-images of the double points between
$M_{ki}:=F_k^{-1}(H_i(D^{2k+1}))$ and $M_{kj}:=F_k^{-1}(H_j(D^{2k+1}))$,
which by transversality may be assumed to be manifolds.
Since $h_0$ is a $(123)$-null-homotopy, $\partial M_{ki}\i S^{2k-1}\x\{1\}$,
and since $h_1$ is a $(321)$-null-homotopy,
$\partial M_{kj}\i S^{2k-1}\x\{0\}$.
Pushing $M_{ki}$ upwards and $M_{kj}$ downwards in $S^{2k-1}\x I$ until they
become disjoint yields an oriented bordism between the double points in
question. \qed
\enddemo

It also follows from Lemma 2.1 and the proof of Theorem 2.2 that $\mu(f)$
can be evaluated using any $(123)$-null-homotopy or $(321)$-null-homotopy
that $f$ admits.

\definition{From ornaments to link maps} An ornament
$f\:S^{2k-1}\sqcup S^{2k-1}\sqcup S^{2k-1}\to\R^{3k-1}$ gives rise to a link
map $\breve f\:S^{2k-1}\sqcup S^{2k-1}\sqcup S^{2k-1}\emb\R^{3k}$ by lifting
the first sphere to overpass the second, the second to overpass the third,
and the third to overpass the first wherever they intersected in $\R^{3k-1}$.
Now a $(123)$-null-homotopy of $\breve f$ can be found within the track
(in the upper half-space of $\R^{3k}$) of any generic homotopy between $f$
and a trivial link map $t$, constant on each component.
So $\breve\mu(f)=\mu(\breve f)$.
This formula is well-known in the classical dimension \cite{Fe}.

Another way to see that $\breve\mu(f)=\mu(\breve f)$ is by noticing that
every $1=2=3$ point in a generic homotopy between $f$ and $f_0$ corresponds
under \ $\breve{ }$ \  to connected summation with the Borromean rings.
Indeed, \ $\breve{ }$ \  when applied to the Borromean ornament $b$ yields
the (slightly smoothed) Borromean rings $B$ --- if instead of the vertical
decomposition $\R^{3k}=\R^{3k-1}\x\R$ we think of the radial decomposition
$\R^{3k}\but\{0\}=S^{3k-1}\x\R$.
\enddefinition

\definition{Reading off the $\mu$-invariant}
There are two ways to generalize the above definition of $\mu(f)$, without
resorting to Lemma 2.1, to arbitrary link maps
$f\:S^{2k-1}\sqcup S^{2k-1}\sqcup S^{2k-1}\to\R^{3k}$ when $k>1$, and
arbitrary h-Brunnian links when $k=1$.

In the first approach, $\mu(f)$ can be defined as above, except that each
null-homotopy is replaced by a null-homology (i.e\. a map of a punctured
$2k$-pseudo\-manifold).
Such null-homologies exist, since each $H^{2k-1}(\R^{3k}\but f(S^{2k-1}))=0$
because the double point locus of $f$ has dimension $\le k-2$ and so
$H_k(f(S^{2k-1}))=0$ unless $k=1$.%
\footnote{Invariance under link homotopy can also be proved in this fashion.
If $F$ is a link homotopy, its double point locus $D$ is
a $(k-1)$-pseudomanifold, so that
$H^{2k}(\R^{3k}\x I\but F(S^{2k-1}\x I))\simeq H_k(F(S^{2k-1}\x I))$ may
be nonzero.
However its generators can be represented by ``characteristic tori''
$S^k\x S^k$ in a neighborhood of $D$, see \cite{HK; \S2}.
So, although we cannot avoid inadmissible intersections between
one immersed cylinder $F(S^{2k-1}\x I)$ and a null-homology for another, they
can be confined to a neighborhood of $D$, which is disjoint from
the null-homology of the third cylinder by general position.}
This definition generalizes the well-known formula \cite{Co; \S5} for the
$\mu$-invariant of a classical h-Brunnian link as the algebraic number of
$1=2=3$ points for Seifert surfaces, bounded by each component in
the complement to the other two.
This formula can be proved for h-Brunnian link maps $3S^{2k-1}\to\R^{3k}$
with arbitrary $k$ along the lines of \cite{Ma1; \S4}.

In the other approach (\cite{Ko6}, \cite{MM}), one considers an arbitrary
triple of null-homotopies of the components, and adds three ``correction
terms'' to the algebraic number of $1=2=3$ points using that the preimage
of each null-homotopy in every other component is null-bordant, since
the two-component sublink is link homotopically trivial.
An advantage of this definition is a straightforward generalization to all
dimensions and a straightforward proof of invariance under link homotopy
\cite{Ko6}.
(The above proof of Theorem 2.2 is a special case of Koschorke's argument;
yet it has the advantage of not breaking the symmetry, as it has already
been broken in our definition.)
It leads to the same $\mu$-invariant as defined in \cite{Ko1}, \cite{Ko7}
for h-Brunnian link maps \cite{Ko6}, and, as is easily seen, to the same
$\mu$-invariant as defined in the previous paragraph.
\enddefinition

\definition{Inducing $\mu$ from a cohomology class of $\Z/3$}
We want to find a Gauss style description for the $\mu$-invariant of
h-Brunnian link maps.
Here is the basic idea.
From the viewpoint of configuration spaces, we should not think of
the complement of a component in $\R^m$ --- we may only think of two
components being disjoint.
Hence each null-homotopy from the initial definition of $\mu$ (for Brunnian
links) becomes a homotopy of the entire link, possibly moving all components.
Now in order to catch the $1=2=3$ point between the null-homotopies as
a generic event in their eyes, we have to treat each homotopy as living in
its own time, and combine them together into a $3$-parametric family
(just like in \cite{Ko6; (4.3)}, albeit for a different purpose).
Thus, we need to rewrite the $(123)$-condition in terms of the new data
to start with.

Let $X=X_1\sqcup X_2\sqcup X_3$ be a polyhedron.
Consider a {\it $\partial I^3$-homotopy} (parametrized by the boundary of
a cube) $h_T\:X\to\R^m$, $T\in\partial I^3$.
We shall call it {\it triangular} if $h_T(X_i\cap X_j)=\emptyset$
whenever $T$ lies on either of the faces $t_j=0$, $t_i=1$ of the cube
$I^3=\{0\le t_1,t_2,t_3\le 1\}$, where $(ijk)=(123)$.
\enddefinition

\proclaim{Lemma 2.3}
If $f\:S^{2k-1}\sqcup S^{2k-1}\sqcup S^{2k-1}\to\R^{3k}$ is an h-Brunnian
link map, there exists a generic cube of homotopies $H_T$, $T\in I^3$,
such that $H_{0,0,0}=f$, $H_{1,1,1}$ is a trivial link map, and the boundary
$\partial I^3$-homotopy is triangular.
Moreover, the algebraic number of $1=2=3$ points in the $I^3$-homotopy equals
$\mu(f)$.
\endproclaim

\demo{Proof} By Lemma 2.1, there exists a link homotopy $h_t$ from $f$ to
a link $l$ whose components bound an admissible triple of null-homotopies.
Consider this triple as a homotopy $g_t$ from $l$ to a trivial link map $t$,
sending each component to a generic point in $\R^{3k}$.
Define $G_{t_1,t_2,t_3}$ by $g_{t_i}$ on the $i$th component $S^{2k-1}$.
This gives the required $I^3$-homotopy except that $G_{0,0,0}=l$ instead
of $f$.
So it remains to modify $G$ near $(0,0,0)$.
To this end let us first extend it via $h_t$ to an $I\vee I^3$-homotopy,
where the whisker $I$ is attached to $I^3$ at $(0,0,0)$.
Now set $H_T=G_{\phi(T)}$, where $\phi\:I^3\to I\vee I^3$ is some
continuous surjection whose composition with the collapse of the whisker
$I\vee I^3\to I^3$ sends each face of the cube to itself. \qed
\enddemo

\bigskip
\centerline{\epsffile{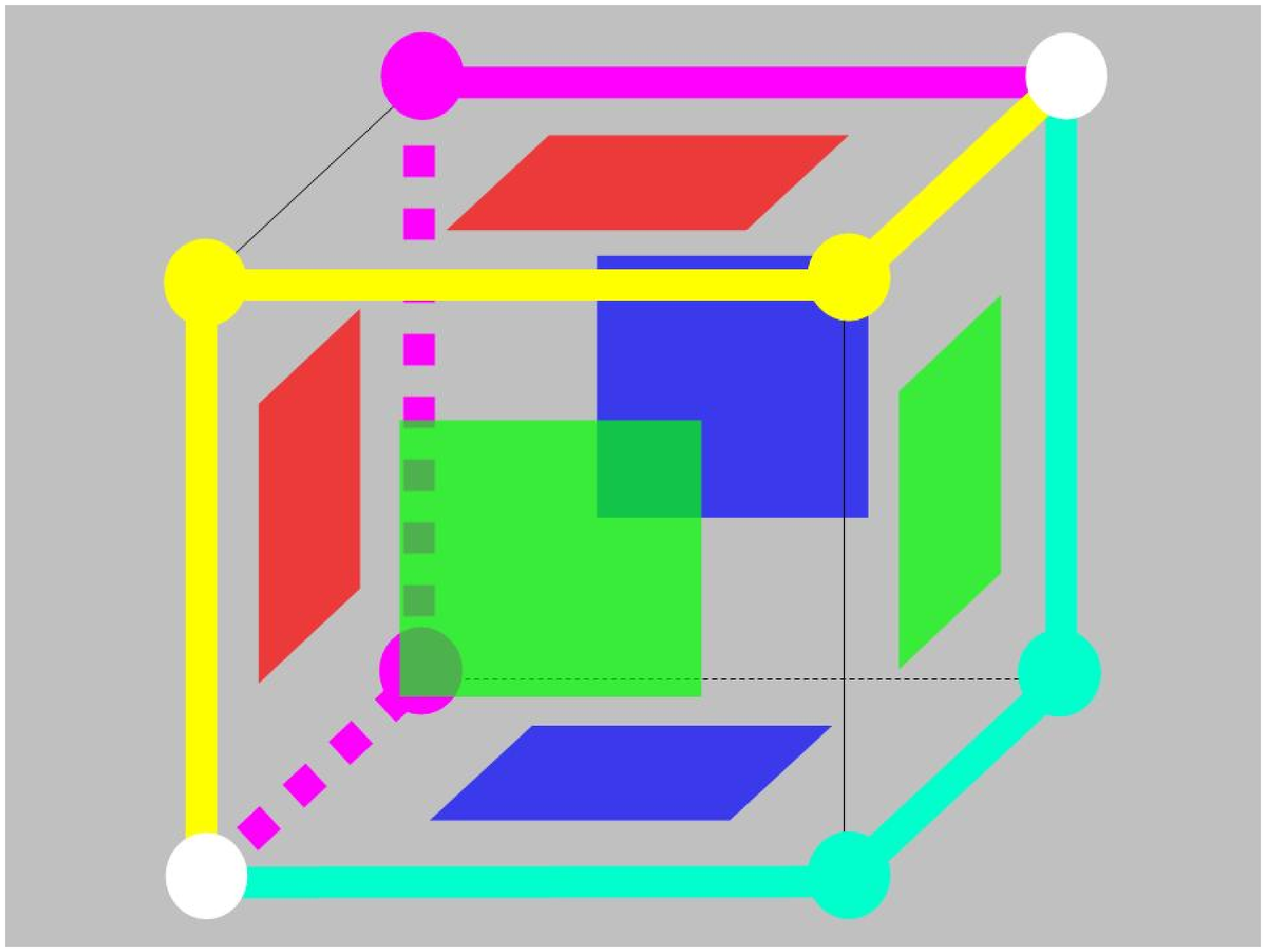}}
\medskip
\centerline{\bf Figure 1: Triangular $\partial I^3$-homotopy}
\bigskip

We shall call a $\partial I^3$-homotopy $h_T\:Q\to(\R^m)^3$ {\it hexagonal}
if $h_T(Q)\cap\Delta^{i=j}=\emptyset$ whenever $T$ lies on either of
the faces $t_k=0$, $t_k=1$ of the cube $I^3=\{0\le t_1,t_2,t_3\le 1\}$
whenever $\{i,j,k\}=\{1,2,3\}$, and additionally $h_T(Q)\i\tl\R^m\Y$
whenever $T$ lies on either edge of the form $t_i=t_j=0$ or $t_i=t_j=1$.
Note that if $h_T$ is hexagonal, its image lies in $\tl\R^m\y$.

\proclaim{Lemma 2.4} Let $G\:Q\to\tl\R^m\Y$ be such that its composition
with each projection $p_i\:\tl\R^m\Y\to\tl\R^m$ forgetting the $i$th point
is null-homotopic.
Then there exists a hexagonal $\partial I^3$-homotopy
$H_T\:Q\to\tl\R^m\y$ such that $H_{0,0,0}=G$ and $H_{1,1,1}$ is a constant
map into $\tl\R^m\Y$.
\endproclaim

\demo{Proof} We first construct a null-homotopy
$H_{T_{ij}(t)}\:Q\to(\R^m)^3\but(\Delta^{k=j}\cup\Delta^{k=i})$ of $G$
corresponding to each of the six shortest paths $T_{ij}$ from $(0,0,0)$
to $(1,1,1)$ in the $1$-skeleton of the cube.

{\sl First edge.}
The fiber $F_i:=p_i^{-1}(x_j,x_k)$ is homeomorphic to $\R^m\but\{x_j,x_k\}$;
we may assume that $x_j-x_k$ is a fixed vector $\vec c\in\R^m$ with
$|\vec c|>1$.
Since $p_iG$ is null-homotopic, $G$ is homotopic to a map $G_i$ into $F_i$.
Writing $\vec e_i$ for the unit vector in the direction of the $t_i$-axis,
let $H_{t\vec e_i}\:Q\to\tl\R^m\Y$ be this homotopy.

{\sl Second edge.}
Let $\pi\:\tl\R^m\to\R^m\but pt$ be the projection along the diagonal,
$\pi(x,y)=x-y$.
Let us extend $p_i$ to the projection
$p_{ij}\:(\R^m)^3\but(\Delta^{k=j}\cup\Delta^{k=i})\to\tl\R^m$.
The fiber $F_{ij}=p_{ij}^{-1}\pi^{-1}(\vec c)$ of $\pi p_{ij}$ contains $F_i$
and coincides with the image of the section $s_{ji}$ of $p_{ji}$, given by
$s_{ji}(x_i,x_k)=(x_i,\,x_k,\,x_k+\vec c)$.
Hence $G_i=s_{ji}p_{ji}G_i$ is homotopic to $s_{ji}p_{ji}G$ with values in
$(\R^m)^3\but(\Delta^{k=j}\cup\Delta^{k=i})$ via the homotopy
$f^{ji}_t:=s_{ji}p_{ji}H_{(1-t)\vec e_i}$.

Let $s_j\:\tl\R^m\to\tl\R^m\Y$ be the section of $p_j$ given by
$s_j(x_i,x_k)=(x_i,x_k,(|x_i|+|x_k|)\vec c)$.
It is homotopic to $s_{ij}$ by a homotopy through sections of $p_{ji}$, given
by combining in order the two homotopies
$$g^{ji}_t(x_i,x_k)=(x_i,\,x_k,\,x_k+[t(|x_i|+|x_k|)+(1-t)]\vec c),$$
$$h^{ji}_t(x_i,x_k)=(x_i,\,x_k,\,(1-t)x_k+(|x_i|+|x_k|)\vec c).$$
Here $g^{ji}_t$ does not touch $\Delta^{k=j}$ since
$t(|x_i|+|x_k|)+(1-t)\ne 0$ for each $t$, and $h^{ji}_t$ since
$|tx_k|\le |x_i|+|x_k|$ for each $t$, and $|\vec c|>1$ (and both homotopies
do not touch $\Delta^{k=i}$ since $x_i\ne x_k$ in the domain).
We define
$H_{\vec e_i+t\vec e_k}\:Q\to(\R^m)^3\but(\Delta^{k=i}\cup\Delta^{k=j})$
by combining in order the homotopies $f^{ji}_t$ and $g^{ji}_tp_{ji}G$ and
$h^{ji}_tp_{ji}G$.

{\sl Third edge.}
Finally, $s_jp_jG$ is null-homotopic via $s_jk^j_t$, where $k^j_t$ is
the given null-homotopy of $p_jG$.
Let $H_{\vec e_i+\vec e_k+t\vec e_j}\:Q\to\tl\R^m\Y$ be this null-homotopy.

{\sl Late face.}
We first observe that $H_{\vec e_i+t\vec e_k}$ followed by
$H_{\vec e_i+\vec e_k+t\vec e_j}$ is homotopic to $H_{\vec e_i+t\vec e_j}$
followed by $H_{\vec e_i+\vec e_j+t\vec e_k}$ through null-homotopies of
$G_i$ with values in $(\R^m)^3\but\Delta^{k=j}$.
Indeed, the former (resp\. the latter) is homotopic to a null-homotopy of
$G_i$ within $F_i\cup\{x_j\}$ (resp\. within $F_i\cup\{x_k\}$), using
a deformation retraction of $F_{ij}$ (resp\. $F_{ik}$) onto $F_i$.
But $F_i\cup\{x_j,x_k\}$ is a Euclidean space disjoint from $\Delta^{k=j}$.

{\sl Early face.}
It remains to check that $H_{t\vec e_i}$ followed by $H_{\vec e_i+t\vec e_k}$
is homotopic to $H_{t\vec e_k}$ followed by $H_{\vec e_k+t\vec e_i}$ through
homotopies between $G$ and $s_jp_jG$ with values in
$(\R^m)^3\but\Delta^{k=i}$.
Let us extend $p_{ji}$ and $p_{jk}$ to the projection
$P_j\:(\R^m)^3\but\Delta^{k=i}\to\tl\R^m$.
Using the homotopy $\phi^{ji}_t$ between the identity map and
$s_{ji}P_j$, given by a deformation retraction of
$(\R^m)^3\but\Delta^{k=i}$ onto $F_{ij}$, we find that $H_{t\vec e_i}$
followed by $f^{ji}_t$ is homotopic to $\phi^{ji}_t G$ through homotopies
between $G$ and $s_{ji}P_jG$.
Finally, $\phi^{ji}_t$ followed by $g^{ji}_tP_j$ and $h^{ji}_tP_j$ and
$h^{jk}_{1-t}P_j$ and $g^{jk}_{1-t}P_j$ and $\phi^{kj}_{1-t}$ is
a self-homotopy of the identity map.
It is null-homotopic since when post-composed with the homotopy equivalence
$P_j$, it equals identically $P_j$ for each value of $t\in [0,6]$. \qed
\enddemo

\definition{Gauss style description of $\mu(f)$}
Consider a polyhedron $X_1\sqcup X_2\sqcup X_3$ and let $Q=X_1\x X_2\x X_3$.
Given an h-Brunnian link map $f\:X_1\sqcup X_2\sqcup X_3\to\R^m$, Lemma 2.4
yields a map $F\:Q\x\partial I^3\to\tl\R^m\y\simeq S^{2m-1}$ such that
$F_{0,0,0}$ is the restriction of $\tl f\y$ and $F_{1,1,1}$ is constant.
Now $(Q\x S^2)/(Q\x pt)$ is homeomorphic to $(Q\sqcup pt)*S^1$, so we have
precisely a well-defined map $\phi\:(Q\sqcup pt)*S^1\to S^{2m-1}$.
Writing $\xi$ for a generator of $H^{2m-1}(S^{2m-1})$, we denote
the image of $\phi^*(\xi)$ under the suspension isomorphism
$H^{2m-1}((Q\sqcup pt)*S^1)\simeq H^{2m-3}(Q)$ by
$$\mu^*(f)\in H^{2m-3}(X_1\x X_2\x X_3).$$
\enddefinition

\proclaim{Theorem 2.5} $\mu^*(f)$ is well-defined up to link homotopy and
equals $2\mu(f)$ when the latter is defined.
\endproclaim

\demo{Proof} Let $f,g\:X_1\sqcup X_2\sqcup X_3\to\R^m$ be link homotopic
h-Brunnian link maps, and let $F_T,G_T\:X_1\x X_2\x X_3\to\tl\R^m\y$ be
hexagonal $\partial I^3$-homotopies with $F_{0,0,0}=f$, $G_{0,0,0}=g$ and
$F_{1,1,1}=G_{1,1,1}=$ a constant map.
We need to show that $F$ and $G$ are homotopic as maps
$X_1\x X_2\x X_3\x\partial I^3\to\tl\R^m\y$.

Let $R$ be the restriction to $\partial I^3$ of a cyclic permutation of
the coordinate axes in $\R^3$.
Since $R$ is isotopic to the identity, $F_T$ is homotopic to $F_{R(T)}$ with
values in $\tl\R^m\y$.
Let $H_T\:X_1\x X_2\x X_3\x\partial I^3\to(\R^m)^3$ be a generic smooth
homotopy between $F_{R(T)}$ and $G_T$.
Similarly to the proof of Theorem 2.2, and using the polyhedral
Pontryagin--Thom construction \cite{BRS}, we can amend $H_T$ so that
its restriction to $X_1\x X_2\x X_3\x\partial I^3$ has values in
$\tl\R^m\y$.

The equation $\mu^*(f)=2\mu(f)$ now follows from Lemma 2.3. \qed
\enddemo

\bigskip
\centerline{\epsffile{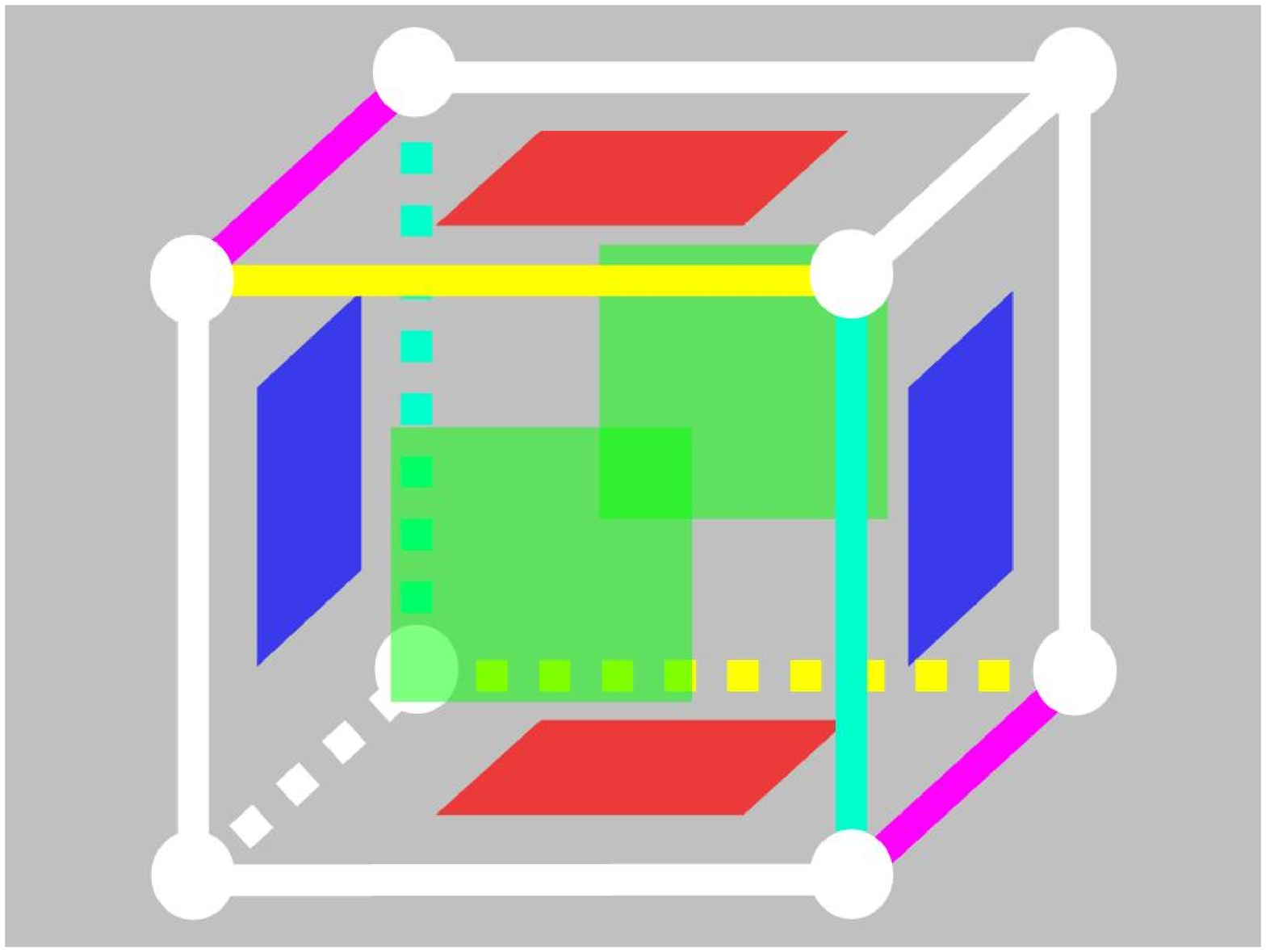}}
\medskip
\centerline{\bf Figure 2: Hexagonal $\partial I^3$-homotopy}
\bigskip

\definition{From link maps to ornaments} The new algebraic definition of
$\mu$ contained essentially in the proof of Lemma 2.4 admits a purely
geometric interpretation.
Given an h-Brunnian a link map $f\:3S^{2k-1}\to\R^{3k}$, we shall construct
geometrically an ornament $J(f)\:3S^{2k+1}\to\R^{3k+2}$ such that
$2\mu(f)=\breve\mu(J(f))$.
Let us put $k=1$ for simplicity of notation.

A circle of null-homotopies of $S^1$ in $S^3$ is a map $S^1*S^1\to S^3*S^1$,
i.e\. $S^3\to S^5$.
Prolonging all null-homotopies slightly in time (by bringing their final
point to the basepoint say), we also get a map $S^3\to S^5\but N(S^1)$,
where $N(S^1)$ is a small regular neighborhood of the distinguished $S^1$.
Thus a triple of circles of null-homotopies for the components of
$f\:S^1\sqcup S^1\sqcup S^1\emb\R^3\subset S^3$ yields a map
$S^3\sqcup S^3\sqcup S^3\to S^5\but N(S^1)\subset\R^5$.
We shall construct just such a triple, or, equivalently, a circle of
homotopies from $f$ to a trivial link map $t$, constant on each component.

A circle of homotopies between $f$ and $t$ is conveniently parametrized
by the boundary of a cube $I^3=\{(t_1,t_2,t_3)\mid 0\le t_i\le 1\}$, with
$f$ corresponding to the vertex $0=(0,0,0)$ and $g$ to
$\vec e_1+\vec e_2+\vec e_3=(1,1,1)$, where $\vec e_i$ are the coordinate
unit vectors in the parameter space.
We shall thus write the homotopy as $H_T$ with $T\in\partial I$.
Let also $\vec v_i$ be the coordinate unit vectors in the ``physical'' space
$\R^3$, and let $B_i$ be the balls of radius $1/3$ centered at
(the endpoints of) $\vec v_i$.
We may assume that the image of $f$ is contained in the ball $B_0$
of radius $1/3$ centered at the origin, and that $t$ sends
the $i$th component to the point $\vec v_i$.
Also let $B_{ij}$, $B_{ijk}$ be the convex hulls of $B_i$ and $B_j$, resp\.
$B_i$, $B_j$ and $B_k$.

{\sl First edge.}
The homotopy along each edge $t\vec v_i$, $t\in I$, is a link homotopy
carrying the $j$th component into $B_j$ and the $k$th component into $B_k$,
where $\{i,j,k\}=\{1,2,3\}$.
Its restriction to the $j$th and $k$th components is provided by
the hypothesis (that $f$ is h-Brunnian) and it can be extended to the
$i$th component by dragging it past the movement of the other components,
which is possible since it never gets close to their double point loci by
general position (this also works in higher dimensions).
We may assume without loss of generality that the homotopy first
``untangles'' the $j$th and $k$th components within $B_0$ and then brings
them into their respective balls, so that the $j$th component always
stays within $B_{0j}$ and the $k$th within $B_{0k}$.
Also without loss of generality, the final map of the homotopy is disjoint
from the points $\vec v_j$ and $\vec v_k$.

{\sl Second edge.}
The homotopy along each edge $\vec v_i+t\vec v_k$ first shrinks the $j$th
component to the point $\vec v_j$ by a radial null-homotopy and then
brings the $i$th and $k$th components back to their original position
via the restriction of $H_{(1-t)\vec v_i}$.
Thus the $j$th component may arbitrarily intersect the $i$th in this
homotopy, but remains disjoint from the $k$th since they stay within
disjoint convex sets $B_j$, $B_{0k}$.
Also the $i$th component is disjoint from the $k$th since they are
immobile in the first stage and move by a link homotopy thereafter.

{\sl Third edge.}
The homotopy along each edge $\vec v_i+\vec v_k+t\vec v_j$ is a link
homotopy keeping the $j$th component fixed (it has been shrunk to $\vec v_j$
already) and shrinking the $i$th to $\vec v_i$ and the $k$th to $\vec v_k$
within $B_{0ik}$.
It is provided by the hypothesis.

{\sl Late face.}
Now the combination of the latter two homotopies is fiberwise linearly
homotopic (through homotopies between $H_{\vec v_i}$ and
$H_{\vec v_i+\vec v_k+\vec v_j}$) to the linear null-homotopy, shrinking
each $l$th component to $v_l$.
Note that the $j$th and $k$th component remain disjoint in this fiberwise
homotopy since they are always contained in disjoint convex sets $B_j$,
$B_{0ik}$.
Along with the symmetric homotopy (with $j$ and $k$ interchanged) this
defines our homotopy along each face $\vec v_i+t\vec v_k+s\vec v_j$.

{\sl Early face.}
Finally, to define it along each face $t\vec v_i+s\vec v_k$,
we similarly note that the combination of $H_{t\vec v_i}$ and
$H_{\vec v_i+t\vec v_k}$ is fiberwise homotopic
(through homotopies between $H_0=f$ and $H_{\vec v_i+\vec v_k}$) to
the homotopy keeping the $i$th and $k$th components fixed and
shrinking the $j$th radially to $\vec v_j$.
Such a fiberwise homotopy can be defined on the $i$th and $k$th
components by $H_{t\vec v_i+s\vec v_k}=H_{(t-s)\vec v_i}$ for $t\ge s$.
Since the $i$th and $k$th components remain disjoint under $H_{t\vec v_i}$,
so they do under this fiberwise homotopy.
The homotopy can be extended to the entire face by defining is symmetrically
(with $i$ and $k$ interchanged) on the other side from the diagonal.
\enddefinition

The procedure $J$ is not entirely canonical, as it slightly depends on
general position arguments.
However, we have

\proclaim{Theorem 2.6} $2\mu(f)=\breve\mu(J(f))$, and $J$ descends to
a well-defined map
$$\CD
\left.\left\{
\text{\rm link maps }S^{2k-1}\sqcup S^{2k-1}\sqcup S^{2k-1}\to\R^{3k}
\right\}\right/
\text{\rm link homotopy}\\
@VVV\\
\left.\left\{
\text{\rm ornaments }S^{2k+1}\sqcup S^{2k+1}\sqcup S^{2k+1}\to\R^{3k+2}
\right\}\right/
\text{\rm ornament homotopy.}
\endCD$$
\endproclaim

\demo{Proof} Since $J(f)$ is essentially a hexagonal
$\partial I^3$-homotopy between $f$ and $t$, the first assertion follows
from Theorem 2.5.
Now the second assertion follows from Theorem 1.1. \qed
\enddemo

Combining $J$ with \ $\breve{ }$ , we get

\proclaim{Corollary 2.7} Let $TLM_k$ be the group of h-Brunnian link maps
$S^{2k-1}\sqcup S^{2k-1}\sqcup S^{2k-1}\to\R^{3k}$ up to link homotopy.
There is a geometrically defined suspension homomorphism
$\sigma\:TLM_k\to TLM_{k+1}$ such that $\mu\sigma=2\mu$.
\endproclaim

Another homomorphism $TLM_k\to TLM_{k+1}$ is a sequence of three
Nezhinskij suspensions.
By \cite{Ko7; Theorem 7.2}, it commutes with the $\mu$-invariant.

It would be interesting to know if $J$ can be decomposed into two operations,
each raising dimension by one.
Note that the $\hat\beta$-invariant from \S3 below vanishes on the $\pm\ne 0$
map given by two null-homotopies of the Borromean rings $B$ along
an appropriate choice of two three-edge paths on the cube.

\head 3. $\beta$-invariant of link maps ($2k-2$ in $3k-2$) \endhead

Given a link map $f\:S^{2k-2}\sqcup S^{2k-2}\to\R^{3k-2}$, one may study
obstructions to embedding a component by a link homotopy.
One way to get such an obstruction is to double the component and consider
an obstruction to removing the intersection between the two copies while
keeping them disjoint from the other original component.

\definition{$\hat\beta$-invariant of $\pm\ne 0$ maps}
Let $f=f_+\sqcup f_-\sqcup f_0\:S^{2k-2}\sqcup S^{2k-2}\sqcup S^{2k-2}
\to\R^{3k-2}$ be a {\it $\pm\ne 0$ map}, i.e\. a map such that the image of
$f_0$ is disjoint from those of $f_+$ and $f_-$.
Let $h_+\sqcup h_-\sqcup h_0\:3D^{2k-1}\to\R^{3k-2}$ be a generic homotopy
from $f$ to a map $t$ sending the components to $3$ distinct points.
We set $\hat\beta(f)$ to be the algebraic number of $1=2=3$ points
between $h_+$, $h_-$ and $f_0$.
This is clearly an invariant of {\it $\pm\ne 0$ homotopy}, i.e.\ homotopy
through $\pm\ne 0$ maps.
\enddefinition

\proclaim{Theorem 3.1} $\hat\beta(f)=\deg(\phi)$, where $\phi$ denotes
the composition
$$S^{2k-2}\x S^{2k-2}\x S^{2k-2}@>f_+\x f_-\x f_0>>(\R^{3k-2})^3\but
(\Delta^{+=0}\cup\Delta^{-=0})@>\simeq>> S^{3k-3}\x S^{3k-3}.$$
\endproclaim

\demo{Proof} We will relate $\hat\beta(f)$ to $\deg(\phi)$ through a sequence
of auxiliary invariants (which might be useful in their own right).
By definition, $\hat\beta(f)$ equals the linking number of $f_0$
and $\partial(h_+\cap h_-)=(f_+\cap h_-)\cup(f_-\cap h_+)$.
This linking number equals the number of $1=2=3$ points between $f_+$, $h_-$
and $h_0$ plus the number of $1=2=3$ points between $f_-$, $h_+$ and $h_0$.
This is the same as $\breve\mu(\hat f)$, where the ornament
$\hat f\:3S^{2k-1}\to\R^{3k-2}\x[-2,2]$ is the track of the following
self-homotopy $\hat f(t)$ of $t=t_+\sqcup t_-\sqcup t_0$:
$$t\quad @<h_-\sqcup h_0<<\quad f_-\sqcup t_+\sqcup f_0\quad @<h_+<<\quad f
\quad @>h_->>\quad f_+\sqcup t_-\sqcup f_0\quad @>h_+\sqcup h_0>>\quad t.$$
Indeed, writing $\hat f=\hat f_+\sqcup\hat f_-\sqcup\hat f_0$, we can extend
$\hat f$ to a generic map $F\:3D^{2k}\to\R^{3k-2}\x[-2,2]\x[0,\infty)$ such
that $F^{-1}(\R^{3k-2}\x[-2,2]\x[1,\infty))$ consists of two $2k$-balls
$D_+$, $D_-$ with
$D_\pm\cap\partial(3D^{2k})=\hat f_\pm^{-1}(\R^{3k-2}\x[\pm1,\pm2])$ and
such that $F^{-1}(\R^{3k-2}\x[0,\pm 2]\x[0,1])$ is sent by $F$ composed with
the projection onto $\R^{3k-2}$ onto the image of
$f_\pm\sqcup h_\mp\sqcup h_0$.
Then the $1=2=3$ points of $F$ project onto those between $f_\pm$, $h_\mp$
and $h_0$.

Now $\breve\mu(\hat f)$ equals the degree of the composition
$$S^{2k-1}\x S^{2k-1}\x S^{2k-1}@>\hat f_+\x\hat f_-\x\hat f_0>>
\tl\R^{3k-1}\y@>\simeq>>S^{6k-3}.$$
This factors through the double suspension of the composition
$$\Sigma(S^{2k-2}\x S^{2k-2}\x S^{2k-2})@>\hat f_+(t)\x\hat f_-(t)\x
\hat f_0(t)>>\tl\R^{3k-2}\y@>\simeq>>S^{6k-5}.$$
Thus $\hat\beta(f)$ equals the degree of the latter composition $\Phi$.
Let us take the homotopies $h_+$ and $h_-$ to be given by translations
(in distinct directions) sufficiently far apart from a cube $I^{3k-2}$
containing the image of $f$, followed by radial null-homotopies.
Then $\Phi$ is the following self-homotopy of a constant map $0$:
$$0@<\ \rho_+\!\!<<\phi_+@<\ \tau_+\!\!<<\phi@>\tau_->>\phi_-@>\rho_->>0,$$
where $\tau_\pm$ is given by a translation of
$(I^{3k-2})^3\but\Delta^{\pm=0}$ within $(\R^{3k-2})^3\but\Delta^{\pm=0}$
along the $\mp$-coordinate, until it is disjoint from
$\Delta^{\mp=0}\cup\Delta^{\mp=\pm}$; and $\rho_\pm$ is given by a
null-homotopy of the translated $(I^{3k-2})^3\but\Delta^{\pm=0}$ within
$(\R^{3k-2})^3\but(\Delta^{\mp=0}\cup\Delta^{\mp=\pm})$.
It follows that $\Phi$ is homotopic to the suspension of $\phi$ composed with
the degree one map $\Sigma(S^{3k-3}\x S^{3k-3})\to S^{6k-5}$ given by
a self-homotopy of a constant map $S^{3k-3}\x S^{3k-3}\to S^{6k-5}$, whose
image is contained in the following subsets of $S^{6k-5}=S^{3k-3}*S^{3k-3}$:
$$*@<D^{3k-3}\x *<<S^{3k-3}\x *@<S^{3k-3}\x D^{3k-2}\!\!\!<<S^{3k-3}\x S^{3k-3}
@>D^{3k-2}\x S^{3k-3}\!\!\!>>*\x S^{3k-3}@>*\x D^{3k-3}\!\!\!>>*.$$
Thus $\deg(\phi)=\deg(\Phi)$. \qed
\enddemo

\definition{$\beta$-invariant of link maps}
Let $g\:S^{2k-2}_*\sqcup S^{2k-2}_0\to\R^{3k-2}$ be a generic link map.
Consider a generic null-homotopy $h\:D^{2k-1}_*\imm\R^{3k-2}$ of
the first component.
If $k$ is even, $\beta(g)\in\Z/2$ is defined to be the parity of the number
of double points between $g(S^{2k-2}_0)$ and the double point manifold
$\Delta(h)$.
If $k$ is odd, the central symmetry of $\R^{3k-2}$ is orientation reversing,
hence the factor exchanging involution of $\R^{3k-2}\x\R^{3k-2}$ reverses
co-orientation of the diagonal, so the factor exchanging involution of
$D^{2k-1}_*\x D^{2k-1}_*$ reverses the co-orientation of
$\hat\Delta(h):=\{(x,y)\mid h(x)=h(y),\, x\neq y\}$ and therefore preserves
its orientation.
Thus $\Delta(h)$, which is the quotient of $\hat\Delta(h)$ by the involution,
is orientable, and we may set $\beta(g)\in\Z$ to be the algebraic number
of double points between $g(S^{2k-2}_0)$ and the double point manifold
$\Delta(h)$.
This is a special case of Koschorke's definition, see \cite{Ko2; proof of
Theorem 4.8}.
(We note that this definition of $\beta$ also applies to link maps of any two
oriented $(2k-2)$-manifolds in $\R^{3k-2}$; the proof of invariance is
straightforward.)
\enddefinition

\definition{$\beta$ and $\hat\beta$} Consider now the projection
$\pi\:S^{2k-2}_+\sqcup S^{2k-2}_-\to S^{2k-2}_*$ and let $f=g\pi$.
If $k$ is odd, $\hat\beta(f)$ coincides by definition with $2\beta(g)$.
Unfortunately, $\beta$ is identically zero for odd $k$ --- in fact,
it is only nonzero for $k=2,4,8$ \cite{Ko2; Theorem 4.8}.
(The extension of $\beta$ to manifold link maps may be nonzero and is
similarly related with the extension of $\hat\beta$ for odd $k$.)
If $k$ is even, $\hat\beta(f)$ is identically zero since it counts every
double point between $g(S^{2k-2}_0)$ and $\Delta(h)$ twice --- with
opposite signs.
\enddefinition

\definition{Gauss style description of $\beta$}
Let $X_1\sqcup X_2$ be a polyhedron and let
$f=f_1\sqcup f_2\:X_1\sqcup X_2\to\R^{3k-2}$ be a link map.
The composition
$$\phi_f\:X_1\x X_2\x X_2@>f_1\x f_2\x f_2>>(\R^{3k-2})^3\but
(\Delta^{1=2}\cup\Delta^{1=3})@>\simeq>> S^{3k-3}\x S^{3k-3}$$
is equivariant with respect to the involution exchanging the second and third
factors of $X_1\x X_2\x X_2$ and the factor exchanging involution of
$S^{3k-3}\x S^{3k-3}$.
The latter is orientation-preserving iff $k$ is odd, so
$H^{6k-6}_{\Z/2}(S^{3k-3}\x S^{3k-3};\,\Z\T^{\otimes (k-1)})\simeq\Z$.
Let $\xi^k$ be a generator of this group, and let us denote $\phi_f^*(\xi^k)$
by $$\beta^*(f)\in H^{6k-6}_{\Z/2}(X_1\x X_2\x X_2;\,\Z\T^{\otimes (k-1)}).$$
\enddefinition

\proclaim{Theorem 3.2} If $X_1=X_2=S^{2k-2}$, then $\beta^*(f)=\beta(f)$.
\endproclaim

\demo{Proof} The factor exchanging involution of $S^{2k-2}\x S^{2k-2}$ is
orientation-preserving, so $H^{6k-6}_{\Z/2}(S^{2k-2}\x S^{2k-2}\x S^{2k-2};\,
\Z\T^{\otimes (k-1)})\simeq\Z$ when $k$ is odd and $\Z/2$ when $k$ is even.
An examination of the proof of Theorem 3.1 shows that it carries over to
the equivariant setting, thus establishing the assertion. \qed
\enddemo

\proclaim{Corollary 3.3} $\beta(f)=0$ if
$f\:S^{2k-2}_*\sqcup S^{2k-2}_0\to\R^{3k-2}$ PL embeds $S^{2k-2}_0$.
\endproclaim

If $f$ PL embeds $S^{2k-2}_*$, this embedding is PL isotopic to
the boundary of an embedded $D^{2k-1}_*$, and it follows from the definition
that $\beta(f)=0$.

A standard proof of Corollary 3.3 is that $\alpha$ is symmetric with respect
to interchange of the components, hence so is $\beta=h(\alpha)$.
The point of the following proof is that it does not use homotopy theory.

\demo{Proof} The embedding of $S^{2k-2}_0$ is PL isotopic to the standard
embedding.
Hence $f$ is link homotopic to the composition of some map
$S^{2k-2}_*\to S^{k-1}$ and the Hopf link
$S^{k-1}\sqcup S^{2k-1}_0\emb\R^{3k-2}$.
Hence the map
$\phi_f\:S^{2k-2}_*\x S^{2k-2}_*\x S^{2k-2}_0\to S^{3k-3}\x S^{3k-3}$
factors up to homotopy through $S^{k-1}\x S^{k-1}\x S^{2k-2}_0$ and so has
zero equivariant degree. \qed
\enddemo

\definition{Sphere instead of the torus}
Let $\vec v$ be a nonzero vector in $\R^m$ and let $\Pi^{\vec v}_2$,
resp\. $\Pi^{\vec v_i}_3$ be the subsets of $(\R^m)^3$ consisting of
all triples $(x,x+\alpha\vec v_i,x)$, resp\. $(x,x,x+\alpha\vec v_i)$
with $\alpha\ge 0$.
Then $(\R^m)^3\but (\Pi^{\vec v}_2\cup\Pi^{\vec v}_3)$ collapses onto
$S^{2m-2}$ equivariantly with respect to the involution exchanging
the second and third factors of $(\R^m)^3$ and the join of the identity
involution on $S^{m-1}$ and the antipodal involution on $S^{2m-3}$.
Moreover, the preimage of $S^{m-1}$ under the collapse is $\Delta^{2=3}$,
and the collapse factors up to equivariant homotopy through the collapse of
$(\R^m)^3\but(\Delta^{1=2}\cup\Delta^{1=3})$ onto $S^{m-1}\x S^{m-1}$
and a degree one map $S^{m-1}\x S^{m-1}\to S^{2m-2}$.
(Note that $\Delta^{2=3}$ cannot be the preimage of the diagonal of
the torus under the latter collapse.)
\enddefinition


\Refs\widestnumber\key{BCSS}


\ref \key Ah \by P. M. Akhmet'ev
\paper On isotopic and discrete realizations of mappings of $n$-sphere to
Euclidean space \jour Mat. Sbornik \vol 187:7 \yr 1996 \pages 3--34
\transl Engl. transl. \jour Sb. Math. \vol 187 \pages 951--980
\endref

\ref \key BRS \by S. Buoncristiano, C. P. Rourke, B. J. Sanderson
\book A Geometric Approach to Homology Theory
\bookinfo London Math. Soc. LNS \vol 18
\publ Cambridge Univ. Press \yr 1976
\endref

\ref \key BT \by A. Bartels, P. Teichner
\paper All two-dimensional links are null-homotopic
\jour Geometry and Topology \vol 3 \yr 1999 \pages 235--252
\endref

\ref \key BCSS \by R. Budney, J. Conant, K. P. Scannell, D. Sinha
\paper New perspectives on self-linking
\jour Adv. Math. \vol 191 \yr 2005 \pages 78--113
\moreref {\tt arXiv:\,math.GT/0303034}
\endref

\ref \key Co \by T. Cochran
\paper Concordance invariance of coefficients of Conway's link polynomial
\jour Invent. \linebreak Math. \vol 82 \yr 1985 \pages 527--541
\endref

\ref \key Fe \by R. Fenn
\book Techniques of Geometric Topology
\bookinfo London Math. Soc. Lect. Note Ser. \vol 57 \yr 1983
\endref

\ref \key FT \by R. Fenn, P. Taylor
\paper Introducing doodles
\inbook Topology of Low-Dimensional Manifolds
\bookinfo Lect. Notes in Math. \vol 722 \yr 1977 \pages 37--43
\endref

\ref \key FM \by W. Fulton and R. MacPherson
\paper A compactification of configuration spaces
\jour Ann. of Math. \vol 139 \yr 1994 \pages 183--225
\endref

\ref \key GKW \by T. G. Goodwillie, J. R. Klein, M. S. Weiss
\paper A Haefliger style description of the embedding calculus tower.
\jour Topology \vol 42 \yr 2003 \pages 509--524
\endref

\ref \key HK1 \by N. Habegger, U. Kaiser
\paper Homotopy classes of $2$ disjoint $2p$-spheres in $\R^{3p+1}$, $p>1$
\jour Topol. Appl. \vol 71 \yr 1996 \pages 1--8.
\endref

\ref \key Hae \by A. Haefliger
\paper Differentiable links
\jour Topology \vol 1 \yr 1962 \pages 241--244
\endref

\ref \key Kh \by M. Khovanov
\paper Doodle groups
\jour Trans. Amer. Math. Soc. \vol 349 \yr 1997 \pages 2297--2315
\endref

\ref \key Ko1 \by U. Koschorke
\paper Higher order homotopy invariants of higher dimensional link maps
\inbook Lecture Notes in Math. \vol 1172 \yr 1985 \pages 116--129
\endref

\ref \key Ko2 \bysame
\paper Link maps and the geometry of their invariants
\jour Manuscr. Math. \vol 61 \yr 1988 \pages 383--415
\endref

\ref \key Ko3 \bysame
\paper On link maps and their homotopy classification
\jour Math. Ann. \vol 286 \yr 1990 \pages 753--782
\endref

\ref \key Ko4 \bysame
\paper Semi-contractible link maps and their suspensions
\inbook Lecture Notes in Math. \vol 1474 \pages 150--169 \yr 1990
\endref

\ref \key Ko5 \bysame
\paper Link homotopy with many components
\jour Topology \vol 30 \yr 1991 \pages 267--281
\endref

\ref \key Ko6 \bysame
\paper Homotopy, concordance and bordism of link maps
\inbook Global analysis in modern mathematics \publ Publish or Perish, Inc.
\publaddr Houston, TX \yr 1993 \pages 283--299
\endref

\ref \key Ko7 \bysame
\paper A generalization of Milnor's $\mu$-invariants to higher-dimensional
link maps
\jour Top\-ology \vol 36 \yr 1997 \pages 301--324
\endref

\ref \key Ma1 \by W. S. Massey
\paper Higher order linking numbers
\inbook Conf. on Algebraic Topology \publ Univ. of Illinois at Chicago Circle
\publaddr Chicago, IL \yr 1969 \pages 174--205
\moreref Reprinted \jour J. Knot Theory Ram. \vol 7 \yr 1998 \pages 393--414
\endref

\ref \key Ma2 \bysame
\paper The homotopy type of certain configuration spaces
\jour Bol. Soc. Mat. Mexicana \vol 37 \yr 1992
\endref

\ref \key M1 \by S. A. Melikhov
\paper Pseudo-homotopy implies homotopy for
singular links of codimension $\ge 3$ \jour Uspekhi Mat. Nauk \vol 55:3
\yr 2000 \pages 183--184
\transl Engl. transl. \jour Russ. Math. Surv. \vol 55 \pages 589--590
\moreref Unabridged version \paper Singular link concordance implies
link homotopy in codimension $\ge 3$ \miscnote preprint, 1998
\endref

\ref \key M2 \by S. A. Melikhov
\paper On maps with unstable singularities
\jour Topol. Appl. \vol 120 \yr 2002 \pages 105--156
\moreref {\tt arXiv:\,math.GT/0101047}
\endref

\ref \key M3 \bysame
\paper Colored finite type invariants and a multi-variable analogue of
the Conway polynomial
\jour {\tt arXiv:\,math.GT/0312007}
\endref

\ref \key M4 \bysame
\paper A compactification of configuration spaces and resolution of
the Thom--Boardman singularities
\miscnote in preparation
\endref

\ref \key M5 \bysame
\paper The van Kampen obstruction and its relatives
\miscnote preprint
\endref

\ref \key M6 \bysame
\paper $n$-Quasi-concordance
\miscnote preprint
\endref

\ref \key MR \bysame, D. Repov\v{s}
\paper $k$-quasi-isotopy: I. Questions of nilpotence
\jour J. Knot Theory Ram. \vol 14 \yr 2005 \pages 571--602
\moreref {\tt arXiv:\,math.GT/0103113}
\endref

\ref \key MS \by S. A. Melikhov, E. V. Shchepin
\paper The telescope approach to embeddability of compacta
\jour {\tt arXiv:\,math.GT/0612085}
\endref

\ref \key MM \by B. Mellor, P. Melvin
\paper A geometric interpretation of Milnor's triple linking numbers
\jour Alg. Geom. Topol. \vol 3 \yr 2003 \pages 557--568
\endref

\ref \key Mer1 \by A. B. Merkov
\paper Finite order invariants of ornaments
\jour J. Math. Sci. \vol 90 \yr 1998 \pages 2215--2273.
\endref

\ref \key Mer2 \bysame
\paper Vassiliev invariants classify plane curves and doodles
\jour Mat. Sbornik \vol 194:9 \yr 2003 \pages 31--62
\transl English transl. \vol 194 \yr 2003 \pages 1301--1333
\endref


\ref \key Si \by D. Sinha
\paper The topology of spaces of knots
\jour {\tt arXiv:\,math.AT/0202287}
\endref

\ref \key Sk$_{\text A}$ \by A. B. Skopenkov
\paper Classification of embeddings outside the metastable dimension
\jour preprint
\endref

\ref \key Sk$_{\text M}$ \by M. B. Skopenkov
\paper A formula for the group of links in the $2$-metastable dimension
\jour preprint
\endref

\ref \key Va1 \by V. A. Vassiliev
\paper Invariants of ornaments
\inbook Singularities and Bifurcations \ed V. I. Arnold
\bookinfo Advances in Soviet Math. \vol 21 \publ Amer. Math. Soc.
\publaddr Providence, RI \yr 1994 \pages 225--262
\endref

\ref \key Va2 \bysame
\book Topology of Complements to Discriminants
\publ Phasis \publaddr Moscow \yr 1997 \lang in Russian
\endref

\ref \key Vo1 \by I. Voli\' c
\paper Finite type knot invariants and calculus of functors
\jour Compos. Math. \vol 142 \yr 2006 \pages 222-�250
\moreref {\tt arXiv:\,math.AT/0401440}
\endref

\ref \key Vo2 \bysame
\paper A survey of Bott--Taubes integration
\jour J. Knot Theory Ram. \toappear
\moreref {\tt arXiv:} {\tt math.GT/0502295}
\endref

\ref \key U \by A. P. Ulyanov
\paper Polydiagonal compactification of configuration spaces
\jour J. Algebr. Geom. \vol 11 \yr 2002 \pages 129--159
\endref

\endRefs
\enddocument
\end